\documentclass{article}
\usepackage{graphicx}
\usepackage{amsmath}
\usepackage{amsthm}
\usepackage{amsfonts}
\usepackage{amssymb}
\newtheorem*{theorem}{Theorem}
\newtheorem*{corollary}{Corollary}
\newtheorem*{definition}{Definition}
\newtheorem*{example}{Example}
\newtheorem*{lemma}{Lemma}
\newtheorem*{notation}{Notation}
\newtheorem*{proposition}{Proposition}
\newtheorem*{remark}{Remark}

\begin{document}

\title{Nonfiliform characteristically nilpotent Lie algebras\footnote{Research
partially supported by the project PB98-0758.}}
\author{Jos\'{e} Mar\'{\i}a Ancochea Berm\'{u}dez\thanks{corresponding author: e-mail: Jose\_Ancochea@mat.ucm.es}
\and Otto Rutwig Campoamor Stursberg\\Departamento de Geometr\'{\i}a y Topolog\'{\i}a\\Facultad CC. Matem\'{a}ticas U.C.M\\28040 Madrid ( Spain )}
\date{}
\maketitle
\begin{abstract}
In this paper we construct large families of characteristically nilpotent Lie
algebras by considering deformations of the Lie algebra $\frak{g}%
_{(m,m-1)}^{4}$ of type $Q_{n}$, and which arises as a central naturally
graded extension of the filiform Lie algebra $L_{n}$. By studying the graded
cohomology spaces we obtain that the sill algebras associated to the models
$\frak{g}_{(m,m-1)}^{4}$ can be interpreted as nilradicals of solvable,
complete Lie algebras. For extreme cocycles we obtain moreover nilradicals of
rigid laws. By considering supplementary cocycles we construct, for any
dimension $n\geq9$, non-filiform characteristically nilpotent Lie algebras and
show that for certain deformations these deformations are compatible with
central extensions.
\end{abstract}

\section*{Introduction}

\bigskip Characteristically nilpotent Lie algebras were born in the late 50's,
as the answer to a question formulated in $1955$ by Jacobson. In [14] he
proved that any Lie algebra defined over a field of characteristic zero and
admitting nondegenerate derivations is nilpotent, and asked for the converse.
Dixmier and Lister gave a negative answer in 1957, constructing an eight
dimensional Lie algebra all whose derivations are nilpotent. The so called
characteristically nilpotent Lie algebras constitute an important subclass
within the nilpotent algebras, and they have been studied in a number of
papers ( [15], [19] )\newline Over the last years the study of
characteristically nilpotent Lie algebras has been reactivated by the fact
that they are of great importance for the topological analysis of the
irreducible components of the variety of nilpotent Lie algebra laws
$\frak{N}^{n}$ [27]. Though there are a lot of constructions of such algebras,
most results correspond to filiform Lie algebras, i. e., algebras with maximal
nilpotence index. In [17] the author studied the graded filiform Lie algebra
$L_{n}$ profoundly, and deduced interesting results about the topology of its
components, as well as giving large families of characteristically nilpotent
Lie algebras.\newline In this paper we study deformations of a family
$\frak{g}_{\left(  m,m-1\right)  }^{4}$ of graded nilpotent Lie algebras of
chracateristic sequence $\left(  n-2,1,1\right)  $, and which are $k$-abelian
of maximal index $k=\left[  \frac{n-2}{2}\right]  $. These algebras are the
simplest case of a wide class of graded Lie algebras called ''of type $Q_{n}%
$'', and which have been introduced in [4]. The family we are considering is
of special interest, as it is the only algebra of type $Q_{n}$ which is a
central, naturally graded extension of the filiform Lie algebra $L_{n}$. The
characteristically nilpotent algebras obtained by deformations of this model
are of value for the general theory, because it is one of the few ( if any )
explicit methods to obtain non-filiform characteristically nilpotent Lie
algebras ( see also [3] ).\newline The paper is divided in four parts. In the
first section we recall the elementary facts about the invariants of nilpotent
Lie algebras used, as well as the rudiments of characteristically nilpotent
Lie algebras. In section two we study the cohomology of nilpotent Lie algebras
and introduce a certain partition of the cohomology groups $H^{2}\left(
\frak{g},\mathbb{C}\right)  $. In the third part we analyze the model
$\frak{g}_{\left(  m,m-1\right)  }^{4}$ and its deformations to obtain, on one
hand, deformations which are isomorphic to the nilradical of solvable,
complete laws, which are rigid for certain cocycles, and on the other,
deformations which are characteristically nilpotent. Finally, we prove that
certain prolongations of these cocycles are compatible with central extensions
of degree one of $\frak{g}_{\left(  m,m-1\right)  }^{4}$ and construct
characteristically nilpotent Lie algebras with characteristic sequence
$\left(  2m-1,2,1\right)  $ for any $m\geq4.$

\section{Generalities}

Let $\frak{L}^{n}$ be the set of complex  Lie algebra laws in dimension $n$.
We identify each law with its structure constants $C_{ij}^{k}$ on a fixed
basis $\left\{  X_{i}\right\}  $ of $\mathbb{C}^{n}$. The Jaciobi identities%
\[
\sum_{l=1}^{n}C_{ij}^{k}C_{kl}^{s}+C_{jk}^{l}C_{il}^{s}+C_{ki}^{l}C_{jl}^{s}=0
\]
for $1\leq i\leq j<k\leq n,\;1\leq s\leq n$ show that $\frak{L}^{n}$ is an
algebraic variety. The nilpotent Lie algebra laws $\frak{N}^{n}$ are a closed
subset in $\frak{L}^{n}$. The linear group $GL\left(  n,\mathbb{C}\right)  $
acts on $\frak{L}^{n}$ by changes of basis. If $\mathcal{O}\left(  \mu\right)
$ denotes the orbit of the law $\mu$ by this action, it is easy to see that it
is a regular subvariety [12].

\begin{definition}
If the orbit $\mathcal{O}\left(  \mu\right)  $ is open in $\frak{L}^{n}$
$\left(  \text{resp. }\frak{N}^{n}\right)  $, the law $\mu$ is rigid in
$\frak{L}^{n}$ $\left(  \text{resp. }\frak{N}^{n}\right)  $.
\end{definition}

For later use, it is convenient to recall some notation :
\[
N_{p}^{n}=\left\{  \mu\in\frak{N}^{n}\;|\;nilindex\left(  \mu\right)  \leq
p-1\right\}
\]

\subsection{\bigskip}

Consider a complex nilpotent Lie algebra $\frak{g}=\left(  \mathbb{C}^{n}%
,\mu\right)  $. For each $X\in\mathbb{C}^{n}$ we denote $c\left(  X\right)  $
the ordered sequence of dimensions of Jordan blocks of the adjoint operator
$ad_{\mu}\left(  X\right)  $.

\begin{definition}
The characteristic sequence of $\frak{g}$ is an isomorphism invariant
$c\left(  \frak{g}\right)  $ defined as
\[
c\left(  \frak{g}\right)  =\sup_{X\in\frak{g}-C^{1}\frak{g}}\left\{  c\left(
X\right)  \right\}
\]
where $C^{1}\frak{g}$ denotes the derived subalgebra.
\end{definition}

\begin{definition}
An $n$-dimensional nilpotent Lie algebra $\frak{g}$\ \ is called filiform if
$c\left(  \frak{g}\right)  =\left(  n-1,1\right)  $.
\end{definition}

Recall that the filiform model algebras $L_{n}$ and $Q_{n}$ are :

\begin{enumerate}
\item  $L_{n}$ is the $\left(  n+1\right)  $-dimensional Lie algebra $\left(  n\geq3\right)$
defined by
\[
[X_{1},X_{i}]=X_{i+1},\;\;2\leq i\leq n
\]
where $\left\{  X_{1},..,X_{n+1}\right\}  $ is a basis.

\item  $Q_{2m-1}$ is the $2m$-dimensional Lie algebra $\left(  m\geq3\right)$ defined by
\[
\begin{array}
[c]{cc}%
\lbrack X_{1},X_{i}]=X_{i+1}, & 2\leq i\leq2m-1\\
\;[X_{j},X_{2m+1-j}]=\left(  -1\right)  ^{j}X_{2m}, & 2\leq j\leq m
\end{array}
\]
where $\left\{  X_{1},..,X_{2m}\right\}  $ is a basis.
\end{enumerate}

\subsection{\bigskip}

\begin{definition}
Let $\frak{g}$ be a solvable, non nilpotent Lie algebra. Then $\frak{g}$ is
called decomposable if
\[
\frak{g}=\frak{n}\oplus\frak{t}%
\]
where $\frak{n}$ is the nilradical of $\frak{g}$ and \ $\frak{t}$ an exterior
torus of derivations, i.e., an abelian subalgebra consisting of $ad$%
-semisimple endomorphisms.
\end{definition}

The structure of these algebras is well known by a theorem of Carles [6] :

\begin{theorem}
A rigid solvable Lie algebra $\frak{g}$ is decomposable.
\end{theorem}

\begin{definition}
If $\frak{t}$ is a maximal torus for the Lie algebra $\frak{g}$, then
$r=\dim\frak{t}$ is called the rank of $\frak{g}$, noted $rank\left(
\frak{g}\right)  $.
\end{definition}

\subsection{\bigskip}

We recall here the elementary facts about characteristically nilpotent Lie algebras.

\begin{definition}
A Lie algebra $\frak{g}$ \ is called characteristically nilpotent if the
algebra of derivations $Der\left(  \frak{g}\right)  $ is nilpotent.
\end{definition}

From this definition it follows that $\frak{g}$ is itself nilpotent.

\begin{proposition}
The characteristically nilpotent Lie algebras constitute a constructible set
in the variety $\frak{N}^{n}$, empty for $n\leq6$ and nonempty for $n\geq7$.
\end{proposition}

A proof of this result can be found in [7].

\section{Cohomology of nilpotent Lie algebras}

Let $\frak{g}$ be a nilpotent Lie algebra whose nilindex is $p$. We have the
ascending and descending central sequences $\left\{  C^{j}\frak{g}\right\}  $
and $\left\{  C_{j}\frak{g}\right\}  $, which are respectively defined by
\[
C^{0}\frak{g=g,\;}C^{q}\left(  \frak{g}\right)  =\left[  C^{q-1}\left(
\frak{g}\right)  ,\frak{g}\right]  ,\;for\;k\geq1
\]%
\[
C_{0}\frak{g}=0,\;C_{q}\left(  \frak{g}\right)  =\left\{  X\in\frak{g\;}%
|\;\left[  X,\frak{g}\right]  \subset C_{q-1}\left(  \frak{g}\right)
\right\}
\]
We set
\[
S_{q}=\frak{g},\;for\;q\leq1,\;\;S_{q}=C^{q-1}\left(  \frak{g}\right)
,\;q\geq1
\]%
\[
T_{q}=\frak{g},\;for\;q\leq1,\;\;T_{q}=C_{p+1-q}\left(  \frak{g}\right)
,\;q\geq1
\]
obtaining two descending filtrations of $\frak{g}$. It is usual to consider
the algebra $\frak{g}$ filtered by $\left\{  S_{q}\right\}  $ and the
$\frak{g}$-module $\frak{g}$ ( respect to the adjoint representation )
filtered by $\left\{  T_{q}\right\}  $. This choice induces descending
filtrations
\[
\left\{  F_{k}Z^{j}\left(  \frak{g},\frak{g}\right)  \right\}  ,\;\left\{
F_{k}H^{j}\left(  \frak{g},\frak{g}\right)  \right\}  ,\;\left\{  F_{k}%
B^{j}\left(  \frak{g},\frak{g}\right)  \right\}
\]
in the cohomology space, which is compatible with the coboundary operator. The
previous filtration has been profoundly studied by Khakimdjanov for the
filiform Lie algebra $L_{n}$ in [15], where he obtained important results on
the topology of the irreducible components of the variety $\frak{N}^{n}$ of
nilpotent Lie algebra laws.

Let
\begin{align*}
d_{i}  &  =\dim_{\mathbb{C}}S_{i},\;1\leq i\leq p\\
d_{p+1}  &  =0
\end{align*}

\begin{lemma}
. \newline 1) If $j\in\mathbb{N}$ and $j>d_{1}$, then $F_{r}Z^{j}\left(
\frak{g},\frak{g}\right)  =Z^{j}\left(  \frak{g},\frak{g}\right)
=0,\;r\in\mathbb{Z}$\newline 2) If $d_{s}<j<d_{s-1}$ for some $1\leq s\leq p$,
then $F_{r}Z^{j}\left(  \frak{g},\frak{g}\right)  =Z^{j}\left(  \frak{g}%
,\frak{g}\right)  $ for $r\leq q$, where
\[
q=-\left[  pd_{p}+\left(  p-1\right)  \left(  d_{p-1}-d_{p}\right)
+..+s\left(  d_{s}-d_{s+1}\right)  +\left(  s-1\right)  \left(  j-1-d_{s}%
\right)  \right]
\]
\end{lemma}

A proof can be found in [16].

\begin{lemma}
\bigskip If $r\leq p\left(  1-j\right)  $ then $F_{r}Z^{j}\left(
\frak{g},\frak{g}\right)  =Z^{j}\left(  \frak{g},\frak{g}\right)  .$
\end{lemma}

\begin{corollary}
If $f$ is a derivation of $\frak{g}$, then $f\left(  S_{i}\right)  \subset
T_{i}$ for any $i.$
\end{corollary}

These results, due originally to Vergne [26], were used by Khakimdjanov [16]
to construct large families of characteristically nilpotent Lie algebras which
are deformations of the filiform model algebra $L_{n}$.

\subsection{\bigskip}

Recall that a central extension of a Lie algebra $\frak{g}$ by $\mathbb{C}%
^{p}$ is an exact sequence of Lie algebras
\[
0\rightarrow\mathbb{C}^{p}\rightarrow\overset{-}{\frak{g}}\rightarrow
\frak{g}\rightarrow0
\]
with $\mathbb{C}^{p}\subset Z\left(  \overset{-}{\frak{g}}\right)  $. If
$\mu_{0}$ is the law of $\frak{g}=\left(  C^{n},\mu_{0}\right)  $, the the law
of $\overset{-}{\frak{g}}$ can be expressed as
\[
\mu\left(  \varphi\right)  \left(  X,Y\right)  =\left\{
\begin{array}
[c]{r}%
\mu_{0}\left(  X,Y\right)  +\varphi\left(  X,Y\right)  ,\;\left(  X,Y\right)
\in\frak{g}^{2}\\
0\;\;,\;\;\left(  X,Y\right)  \in\mathbb{C}^{p}\times\mathbb{C}^{m}%
\end{array}
\right.  \;
\]
where $\varphi\in Z^{2}\left(  \frak{g},\mathbb{C}^{p}\right)  $. The orbit
$\mu\left(  Z^{2}\left(  \frak{g},\mathbb{C}^{p}\right)  \right)  $ by the
action of the group $GL\left(  m+p,\mathbb{C}\right)  $ consists of the law
over $\mathbb{C}^{m+p}$ obtained by this manner.

The following result, due to Carles [6], is of importance in the study of
characteristically nilpotent Lie algebras :

\begin{proposition}
The subspace $E_{c,p}\left(  \frak{g}\right)  $ of laws in $\frak{N}^{m+p}$
which are ( up to isomorphism ) central extensions of $\frak{g}$ by
$\mathbb{C}^{p}$ is constructible, irreducible and of dimension
\[
m\left(  m+p\right)  -\dim Der\left(  \frak{g}\right)  +p\dim Z^{2}\left(
\frak{g},\mathbb{C}^{p}\right)  -\rho
\]
where $\rho$ is the minimum for $\overset{-}{\frak{g}}\in E_{c,p}\left(
\frak{g}\right)  $ of the dimensions of subspaces of the Grasmannian
$Gr_{p}\left(  Z\left(  \overset{-}{\frak{g}}\right)  \right)  $ consisting of
ideals $\frak{a}\subset Z\left(  \overset{-}{\frak{g}}\right)  $ such that
$\frac{\overset{-}{\frak{g}}}{\frak{a}}\simeq\frak{g}$.
\end{proposition}

\begin{proof}
The space $E_{c,p}\left(  \frak{g}\right)  $ is the image of $Gl\left(
m+p,\mathbb{C}\right)  \times Z^{2}\left(  \frak{g},\mathbb{C}^{p}\right)  $
by the algebraic map%
\[
\left(  s,\varphi\right)  \mapsto s\ast\mu\left(  \varphi\right)
\]
where * denotes the action of $Gl\left(  m+p,\mathbb{C}\right)  $ on
$\frak{N}^{m+p}$. It follows that $E_{c,p}\left(  \frak{g}\right)  $ is
irreducible and constructible. If $\phi$ is the law of $\overset{-}{\frak{g}}%
$, then its inverse image is algebraically isomorphic to the first projection
\[
X\left(  \phi\right)  =\left\{  s\in GL\left(  m+p,\mathbb{C}\right)
\;|\;s^{-1}\ast\phi\in\mu\left(  Z^{2}\left(  \frak{g},\mathbb{C}^{p}\right)
\right)  \right\}
\]
The algebraic mapping which associates $s\left(  \mathbb{C}^{p}\right)  $ to
$s\in X\left(  \phi\right)  $ in the Grasmannian $Gr_{p}\left(  \mathbb{C}%
\right)  $ by the image of the set of central ideals of $\overset{-}{\frak{g}%
}$ \ whose factor is isomorphic to $\frak{g}$. If $s_{1}$ and $s_{2}$ give the
same image $\frak{h}$ then $s_{1}^{-1}s_{2}$ belongs to the stabilizer of
$\mathbb{C}^{p}$, mapping $\mu\left(  \varphi_{2}\right)  $ on $\mu\left(
\varphi_{1}\right)  $. We deduce that the minimum of dimensions for $X\left(
\phi\right)  $ equals $\dim Der\left(  \frak{g}\right)  +mp+\rho$.
\end{proof}

As known, the second cohomology space with values in $\mathbb{C}$ is obtained
as follows : Let $\lambda:\bigwedge^{2}\frak{g\rightarrow g}$ be the linear
mapping defined by
\[
\lambda\left(  X\wedge Y\right)  =\mu_{0}\left(  X,Y\right)  ,\;\;X,Y\in
\frak{g}%
\]
and let $\Omega$ be the vector subspace of $\bigwedge^{2}\frak{g}$ generated
by the elements
\[
\left[  X,Y\right]  \wedge Z+\left[  Y,Z\right]  \wedge X+\left[  Z,X\right]
\wedge X,\;\;\;\left(  X,Y,Z\right)  \in\frak{g}^{3}%
\]
It follows that the second homology space $H_{2}\left(  \frak{g}%
,\mathbb{C}\right)  $ coincides with the factor $\frac{Ker\lambda}{\Omega}$,
and by duality%
\[
H^{2}\left(  \frak{g},\mathbb{C}\right)  =Hom\left(  \frac{Ker\lambda}{\Omega
},\mathbb{C}\right)  =H_{2}\left(  \frak{g},\mathbb{C}\right)  ^{\ast}%
\]

This space is of special interest in the determination of central extensions
of degree $1$. Now, let $\left\{  X_{1},..,X_{n}\right\}  $ be a basis for the
Lie algebra $\frak{g}$ and define the cocycles $\varphi_{ij}\in H^{2}\left(
\frak{g},\mathbb{C}\right)  ,\;\left(  \;i,j,k,l\leq n\right)  $
\[
\varphi_{ij}\left(  X_{k},X_{l}\right)  =\left\{
\begin{array}
[c]{cc}%
\delta_{ij}\in\mathbb{C} & \text{if \ }i=k,j=l\\
0 & \text{otherwise}%
\end{array}
\right.
\]
where $\delta_{ij}$ denotes the Kronecker delta function. It is easy to verify
that $\sum_{i,j}a^{ij}\varphi_{ij}=0$ \ $\left(  a_{ij}\in\mathbb{C}\right)  $
if and only if $\sum a^{ij}\left(  X_{i}\wedge X_{j}\right)  \in\Omega$. \newline 

\begin{notation}
For $k\geq2$ let
\[
H_{k}^{2,t}\left(  \frak{g},\mathbb{C}\right)  =\left\{  \varphi_{ij}\in
H^{2}\left(  \frak{g},\mathbb{C}\right)  \;|\;i+j=2t+1+k\right\}  ,\;1\leq
t\leq\left[  \frac{\dim\frak{g}-3}{2}\right]
\]%
\[
H_{k}^{2,\frac{t}{2}}\left(  \frak{g},\mathbb{C}\right)  =\left\{
\varphi_{ij}\in H^{2}\left(  \frak{g},\mathbb{C}\right)
\;|\;i+j=t+1+k\right\}  ,\;1\leq t\leq\left[  \frac{\dim\frak{g}-3}{2}\right]
,t\equiv1\left(  \operatorname{mod}2\right)
\]
\end{notation}

If $\mathbf{E}_{c,1}\left(  \frak{g}\right)  $ defines the subset of
$E_{c,1}\left(  \frak{g}\right)  $ formed by the central extensions of
$\frak{g}$ by $\mathbb{C}$ which are additionally naturally graded, we can
define the sets
\[
\mathbf{E}_{c,1}^{t,k_{1},..,k_{r}}\left(  \frak{g}\right)  =\left\{  \mu\in
E_{c,1}(\frak{g)\;}|\;\mu=\mu_{0}+\left(  \sum a_{kij}\varphi_{ij}^{k}\right)
,\;\varphi_{ij}\in H_{k}^{2,t}\left(  \frak{g},\mathbb{C}\right)  ,a_{kij}%
\in\mathbb{C}\right\}
\]%
\[
\mathbf{E}_{c,1}^{\frac{t}{2},k_{1},..,k_{r}}\left(  \frak{g}\right)
=\left\{  \mu\in E_{c,1}(\frak{g)\;}|\;\mu=\mu_{0}+\left(  \sum a_{kij}%
\varphi_{ij}^{k}\right)  ,\;\varphi_{ij}\in H_{k}^{2,\frac{t}{2}}\left(
\frak{g},\mathbb{C}\right)  ,a_{kij}\in\mathbb{C}\right\}
\]
If $\left\{  X_{1},..,X_{n},X_{n+1}\right\}  $ is a basis of $\overset
{-}{\frak{g}}\in E_{c,1}\left(  \frak{g}\right)  $ the its law $\mu$ is given
by :
\[
\mu\left(  X,Y\right)  :=\mu_{0}\left(  X,Y\right)  +\left(  \sum
a_{kij}\varphi_{ij}^{k}\right)  X_{n},\;\;\left(  X,Y\right)  \in\frak{g}^{2}%
\]

\begin{lemma}
The following set-theoretically identity holds%
\[
\mathbf{E}_{c,1}\left(  \frak{g}\right)  =\bigsqcup_{t,k_{i}}\left(
\mathbf{E}_{c,1}^{t,k_{1},..,k_{r}}\left(  \frak{g}\right)  \cup
\mathbf{E}_{c,1}^{\frac{t}{2},k_{1},..,k_{r}}\left(  \frak{g}\right)
\right)
\]
\end{lemma}

The proof is trivial.

\begin{example}
For the Lie algebra $\frak{g}_{7}$ defined by the brackets%
\begin{align*}
\lbrack X_{1},X_{i}] &  =X_{i+1},\;2\leq i\leq5\\
\lbrack X_{2},X_{3}] &  =X_{7}%
\end{align*}
we have the decomposition
\[
\mathbf{E}_{c,1}\left(  \frak{g}_{7}\right)  =\mathbf{E}_{c,1}^{\frac{5}{2}%
,1}\left(  \frak{g}_{7}\right)  \cup\mathbf{E}_{c,1}^{\frac{3}{2},3}\left(
\frak{g}_{7}\right)  \cup\mathbf{E}_{c,1}^{\frac{3}{2},4}\left(  \frak{g}%
_{7}\right)  \cup\mathbf{E}_{c,1}^{\frac{3}{2},4,5}\left(  \frak{g}%
_{7}\right)
\]
\end{example}

\begin{remark}
Recall that the $2$-cocycles $\psi\in Z^{2}\left(  \frak{g},\frak{g}\right)  $
of the Chevalley cohomology can be interpreted as the infinitesimal
deformations of the Lie algebra $\frak{g}$.
\end{remark}

\begin{definition}
Let\ $\psi\in Z^{2}\left(  \frak{g},\frak{g}\right)  $ be a cocycle of the Lie
algebra $\frak{g}=\left(  \mathbb{C}^{n},\mu\right)  $. Then $\psi$ is called
linearly expandable if the operation
\[
\left(  \mu+\psi\right)  \left(  X,Y\right)  =\mu\left(  X,Y\right)
+\psi\left(  X,Y\right)
\]
defines a Lie algebra law over $\mathbb{C}^{n}$.
\end{definition}

\section{The algebra $\frak{g}_{\left(  m,m-1\right)  }^{4}$}

\bigskip In this section we are interested on a certain class of nilpotent Lie
algebras, which are called ''Lie algebras of type $Q_{n}$''.

\begin{definition}
A nilpotent Lie algebra $\frak{g}$ is called $k$-abelian if the ideal
$C^{k}\frak{g}$ of the central descending sequence is abelian and
$C^{k-1}\frak{g}$ is not abelian. \newline The index $k$ is called abelianity
index of the Lie algebra $\frak{g}$.
\end{definition}

\begin{remark}
If $p$ is the nilindex of $\frak{g}$, it is trivial to verify that
\[
k\leq\left[  \frac{p}{2}\right]
\]
where $\left[  \left(  {}\right)  \right]  $ is the integer part function.
\end{remark}

\begin{definition}
A $k$-abelian nilpotent Lie algebra of nilindex $p$ is called of type $Q_{n}$
if $k=\left[  \frac{p}{2}\right]  $.
\end{definition}

The denomination is inspired on the following property of the filiform algebra
$Q_{n}$ :

\begin{proposition}
Let $\frak{g}$ be a naturally graded, $n$-dimensional filiform Lie algebra.
Then $\frak{g}$ is $k$-abelian of maximal abelianity index if and only
\newline 1) $\dim\frak{g}=2m$\newline 2) $\frak{g}\simeq Q_{2m-1}.$
\end{proposition}

\begin{remark}
Vergne proved in $1966$ that a naturally graded filiform Lie algebra is either
isomorphic to $L_{n}$ or $Q_{n}$, and where the last algebra exists only in
even dimension. As $L_{n}$ is clearly $1$-abelian, a $k$-abelian filiform
algebra with maximal abelianity index must be an isomorphic copy of $Q_{n}$.
This forces the even dimensionality.
\end{remark}

For $m\geq4$ let $\frak{g}_{\left(  m,m-1\right)  }^{4}$ be the Lie algebra
whose structural equations are
\begin{align*}
d\omega_{1}  &  =d\omega_{2}=0\\
d\omega_{j}  &  =\omega_{1}\wedge\omega_{j-1},\;3\leq j\leq2m\\
d\omega_{2m+1}  &  =\sum_{j=2}^{m}\left(  -1\right)  ^{j}\omega_{j}%
\wedge\omega_{2m+1-j}%
\end{align*}
where $\left\{  \omega_{1},..,\omega_{2m+1}\right\}  $ is a basis of $\left(
\mathbb{C}^{2m+1}\right)  ^{\ast}$.

\begin{proposition}
For any $m\geq4$ the Lie algebra $\frak{g}_{\left(  m,m-1\right)  }^{4}$ is
naturally graded, $\left(  m-1\right)  $-abelian of characteristic sequence
$\left(  2m-1,1,1\right)  $.
\end{proposition}

\begin{proof}
We only observe that the $\left(  m-1\right)  $-abelianity is given by the
differential form $d\omega_{2m+1}$..
\end{proof}

\begin{notation}
Recall that
\[
\delta N_{p}^{n}=\left\{  \mu\in N_{p}^{n}\;|\;nilindex\left(  \mu\right)
=p-1\right\}
\]
\end{notation}

\begin{theorem}
For $m\geq4$ any naturally graded, central extension of $L_{2m-1}$ by
$\mathbb{C}$ whose nilindex is $\left(  2m-1\right)  $ is isomorphic to
$\frak{g}_{\left(  m,m-1\right)  }^{4}$.
\end{theorem}

\begin{proof}
We analyze the elements of the set $\mathbf{E}_{c,1}^{t,2}\left(
L_{2m-1}\right)  $ which have the prescribed nilindex.\newline For $t\neq m-1$
it is immediate to see that the nonexistence of such an extension. For $t=m-1$
the cocycles $\varphi_{ij}$ defining the extension have to satisfy the
relation
\[
i+j=2m+1
\]
The only ones which verify this requirement are the cocycles%
\[
\varphi_{j,2m+1-j},\;2\leq j\leq\left[  \frac{2m+1}{2}\right]
\]
Moreover, if $\left\{  X_{1},..,X_{2m+1}\right\}  $ is the dual basis of
$\left\{  \omega_{1},..,\omega_{2m+1}\right\}  $ we have
\[
X_{2}\wedge X_{2m-1}+\left(  -1\right)  ^{j}X_{j}\wedge X_{2m+1-j}\in
\Omega,\;3\leq j\leq\left[  \frac{2m+1}{2}\right]
\]
Observe that the space $\frac{Ker\lambda}{\Omega}$ gives the presentation
\[
\left\langle X_{1},..,X_{2m+1}\;|\;X_{2}\wedge X_{2m-1}+\left(  -1\right)
^{j}X_{j}\wedge X_{2m+1-j}=0,\;3\leq j\leq\left[  \frac{2m+1}{2}\right]
\right\rangle
\]
and the underlying Lie algebra is obviously isomorphic to $\frak{g}_{\left(
m,m-1\right)  }^{4}$ for any $m\geq4$.
\end{proof}

\begin{corollary}
For any values $\left(  t,k_{1},..,k_{r}\right)  \neq\left(
3,2,0,..,0\right)  ,\;$ $\left(  r\geq1\right)  $ $\ $we have
\[
\mathbf{E}_{c,1}^{t,k_{1},..,k_{r}}\left(  \frak{g}_{\left(  m,m-1\right)
}^{4}\right)  =\emptyset
\]
\end{corollary}

\begin{remark}
It is routine to verify that $\frak{g}_{\left(  m,m-1\right)  }^{4}$
satisfies
\[
S_{q}=T_{q}=\bigoplus_{k\geq q}\frak{g}_{k}%
\]
where $\frak{g}_{k}$ denotes the $k^{th}$ block of the associated graduation
of $\frak{g}_{\left(  m,m-1\right)  }^{4}$.
\end{remark}

Let $\left\{  X_{1},..,X_{2m+1}\right\}  $ denote the dual basis of $\left\{
\omega_{1},..,\omega_{2m+1}\right\}  $ given for $\frak{g}_{\left(
m,m-1\right)  }^{4}$.

\begin{proposition}
The linear mappings $ad\left(  X_{1}\right)  ,..,ad\left(  X_{2m-1}\right)
,f_{1}^{1},f_{1}^{2},f_{1}^{2m+1},$\newline $f_{2}^{2},f_{2}^{3+j}\;\left(
1\leq j\leq2m-4\right)  ,f_{2}^{2m},f_{2}^{2m+1},\;$ where%
\[
\left\{
\begin{array}
[c]{c}%
f_{1}^{1}\left(  X_{1}\right)  =X_{1},\;f_{1}^{1}\left(  X_{2}\right)  =0,\\
f_{1}^{1}\left(  X_{i}\right)  =\left(  i-2\right)  X_{i},\;3\leq i\leq2m,\\
\;f_{1}^{1}\left(  X_{2m+1}\right)  =\left(  2m-3\right)  X_{2m+1}%
\end{array}
\right.
\]%
\[
\left\{
\begin{array}
[c]{c}%
f_{2}^{2}\left(  X_{1}\right)  =0\\
\;f_{2}^{2}\left(  X_{i}\right)  =X_{i},\;2\leq i\leq2m,\\
\;f_{2}^{2}\left(  X_{2m+1}\right)  =2X_{2m+1}%
\end{array}
\right.
\]%
\[
\left\{  f_{1}^{2}\left(  X_{1}\right)  =X_{2},\;f_{1}^{2}\left(
X_{2}\right)  =X_{2m+1}\right.
\]%
\[
\left\{  \;\;f_{1}^{2m+1}\left(  X_{1}\right)  =X_{2m+1}\right.
\]%
\[
\left\{  f_{2}^{3+j}\left(  X_{k}\right)  =X_{k+j+1},\;1\leq j\leq2m-4,\;2\leq
k\leq2m-1-j\;\right.
\]%
\[
\left\{  f_{2}^{2m}\left(  X_{2}\right)  =X_{2m};\;\;f_{2}^{2m+1}\left(
X_{2}\right)  =X_{2m+1}\right.
\]
\newline and the undefined images are zero, form a basis of $Der\left(
\frak{g}_{\left(  m,m-1\right)  }^{4}\right)  $.
\end{proposition}

\begin{proof}
The listed endomorphisms are clearly derivations of \ $\frak{g}_{\left(
m,m-1\right)  }^{4}$. Their linear independence is also easy to verify. As
$\frak{g}_{\left(  m,m-1\right)  }^{4}$ is naturally graded, a derivation $d$
decomposes as
\[
d=d_{0}+d_{1}+...+d_{2m-2}%
\]
where $d_{j}\in Der\left(  \frak{g}_{\left(  m,m-1\right)  }^{4}\right)  $ and
$d_{j}\left(  \frak{g}_{k}\right)  \subset\frak{g}_{k+j}$ for $1\leq
j\leq2m-2,\;1\leq k\leq2m-1$. From this decomposition it is not difficult to
prove that
\begin{align*}
d_{0} &  =a_{0}^{1}f_{1}^{1}+a_{0}^{2}f_{2}^{2}+a_{0}^{3}f_{1}^{2}\\
d_{1} &  =a_{1}^{1}ad\left(  X_{1}\right)  +a_{2}^{2}ad\left(  X_{2}\right)
\\
d_{k} &  =a_{k}^{1}ad\left(  X_{k+1}\right)  +a_{k}^{2}f_{2}^{k+2},\;\;2\leq
k\leq2m-3\\
d_{2m-2} &  =a_{2m-2}^{1}ad\left(  X_{2m-1}\right)  +a_{2m-2}^{2}f_{2}%
^{2m}+a_{2m-2}^{3}f_{2}^{2m+1}+a_{2m-2}^{4}f_{1}^{2m+1}%
\end{align*}
where $a_{i}^{j}\in\mathbb{C}$ for $0\leq i\leq2m-2,\;j=1,2,3,4$.
\end{proof}

\begin{corollary}
For $m\geq4,\;\dim Der\left(  \frak{g}_{\left(  m,m-1\right)  }^{4}\right)
=4m+1.$ $\;$
\end{corollary}

\begin{corollary}
For $m\geq4,\;\dim H^{1}\left(  \frak{g}_{\left(  m,m-1\right)  }^{4}%
,\frak{g}_{\left(  m,m-1\right)  }^{4}\right)  =2m.$
\end{corollary}

In particular, we see that for any $m\geq4$ the rank of the Lie algebra
$\frak{g}_{\left(  m,m-1\right)  }^{4}$ is $2$.

Let $H^{2}\left(  \frak{g}_{\left(  m,m-1\right)  }^{4},\frak{g}_{\left(
m,m-1\right)  }^{4}\right)  $ be the second Chevalley cohomology space. We
consider cocicles $\varphi$ which satisfy the following property :

\begin{enumerate}
\item  For any $X\in\frak{g}_{\left(  m,m-1\right)  }^{4}$ such that there
exists $Y\in Z\left(  \frak{g}_{\left(  m,m-1\right)  }^{4}\right)  $ with
$Y\in\frak{g}_{2m-1}$ and $Y\notin im\,ad\left(  X\right)  $ we have
\[
\varphi\left(  X,\frak{g}_{\left(  m,m-1\right)  }^{4}\right)  =0
\]

\item  If $Y\in Z\left(  \frak{g}_{\left(  m,m-1\right)  }^{4}\right)  $ is
such that $Y\in\frak{g}_{2m-1}$ then $Y\notin im\left(  \varphi\right)  $.
\end{enumerate}

We say that the cocycle $\varphi$ satisfies the property $\left(  P\right)  $
and write
\[
\widehat{Z}^{2}\left(  \frak{g}_{\left(  m,m-1\right)  }^{4},\frak{g}_{\left(
m,m-1\right)  }^{4}\right)  =\left\{  \varphi\in Z^{2}\left(  \frak{g}%
_{\left(  m,m-1\right)  }^{4},\frak{g}_{\left(  m,m-1\right)  }^{4}\right)
\;|\;\varphi\text{ satisfies }\left(  P\right)  \right\}
\]
Then, for each $k$, we obtain the subspaces
\[
\left\{  F_{k}\widehat{Z}^{2}\left(  \frak{g}_{\left(  m,m-1\right)  }%
^{4},\frak{g}_{\left(  m,m-1\right)  }^{4}\right)  \right\}  ,\;\left\{
F_{k}\widehat{B}^{2}\left(  \frak{g}_{\left(  m,m-1\right)  }^{4}%
,\frak{g}_{\left(  m,m-1\right)  }^{4}\right)  \right\}  ,\;
\]
and
\[
\left\{  F_{k}\widehat{H}^{2}\left(  \frak{g}_{\left(  m,m-1\right)  }%
^{4},\frak{g}_{\left(  m,m-1\right)  }^{4}\right)  \right\}
\]

For each $m\geq4$ and $0\leq r\leq2m-4,\;$ $2\leq k\leq\left[  \frac{2m-r}%
{2}\right]  $ we define the cocycles
\[
\psi_{k,r}\left(  X_{l},X_{m}\right)  =\left\{
\begin{array}
[c]{cc}%
X_{k+j-1+r} & \text{if \ }l=k,m=j,\;1+k\leq j\leq2m+1-k-r\\
0 & \text{otherwise}%
\end{array}
\right.
\]

\begin{proposition}
For any $m\geq4,$ \ the cocycles \ $\psi_{k,r}$ $\left(  0\leq r\leq
2m-4,\;2\leq k\leq\left[  \frac{2m-r}{2}\right]  \right)  $ form a basis of
$F_{0}\widehat{Z}^{2}\left(  \frak{g}_{\left(  m,m-1\right)  }^{4}%
,\frak{g}_{\left(  m,m-1\right)  }^{4}\right)  $.\newline Moreover, for any
$0\leq j_{0}\leq2m-4,$ \ a basis of $\widehat{Z}_{j_{0}}^{2}\left(
\frak{g}_{\left(  m,m-1\right)  }^{4},\frak{g}_{\left(  m,m-1\right)  }%
^{4}\right)  $ is given by the cocycles $\psi_{k,j_{0}}\;\left(  2\leq
k\leq\left[  \frac{2m-j_{0}}{2}\right]  \right)  $.
\end{proposition}

\begin{proof}
Let $\psi\in F_{0}\widehat{Z}^{2}\left(  \frak{g}_{\left(  m,m-1\right)  }%
^{4},\frak{g}_{\left(  m,m-1\right)  }^{4}\right)  $ be a cocycle. Then, for
$2\leq j\leq m-2,\;0\leq h\leq2m-j$, $1\leq t\leq2m-2j-h$ \ we have
\[
\psi\left(  X_{j},X_{j+1+h}\right)  =\alpha_{j,j+1+t+h}^{2j+t+h}X_{2j+t+h}%
\]
The cocycle condition $\delta\left(  \psi\right)  =0$ implies the relations
\[
\alpha_{j,j+1+h}^{2j+h}=\alpha_{j,j+1+t+h}^{2j+t+h}%
\]
for $2\leq j\leq m-2,\;0\leq h\leq2m-j,\;1\leq t\leq2m-2j-h$, and
\begin{align*}
\alpha_{m-1,m}^{2m-2} &  =\alpha_{m-1,m+1}^{2m-1}\\
\alpha_{m-1,m}^{2m-1} &  =\alpha_{m-1,m+1}^{2m}\\
\alpha_{m,m+1}^{2m} &  =\alpha_{m-1,m}^{2m-2}-\alpha_{m-1,m+2}^{2m}%
\end{align*}
Thus the cocycle $\psi$ can be rewritten as :
\begin{align*}
\psi &  =\sum_{j=2}^{m-2}\left(  \alpha_{j,j+1}^{2j}\psi_{j,0}+...+\alpha
_{j,j+1}^{2m}\psi_{j,2m+1-j}\right)  +\alpha_{m-1,m}^{2m-2}\psi_{m-1,0}%
+\alpha_{m-1,m}^{2m-1}\psi_{m-1,m}+\\
&  +\alpha_{m-1,m}^{2m}\psi_{m-1,2}+\alpha_{m,m+1}^{2m}\psi_{m,0}%
\end{align*}
So the cocycles $\psi_{k,r}\;\left(  0\leq r\leq2m-4,\;2\leq k\leq\left[
\frac{2m-r}{2}\right]  \right)  $ constitute a generator system for
$F_{0}\widehat{Z}^{2}\left(  \frak{g}_{\left(  m,m-1\right)  }^{4}%
,\frak{g}_{\left(  m,m-1\right)  }^{4}\right)  $. Clearly they are linearly
independent, so they form a basis. This proves the first assertion. The second
follows immediately.
\end{proof}

\begin{corollary}
For $m\geq4$ and $0\leq t\leq m-3,$%
\[
\dim\widehat{Z}_{0}^{2}\left(  \frak{g}_{\left(  m,m-1\right)  }^{4}%
,\frak{g}_{\left(  m,m-1\right)  }^{4}\right)  =m-1
\]%
\[
\dim\widehat{Z}_{2t+1}^{2}\left(  \frak{g}_{\left(  m,m-1\right)  }%
^{4},\frak{g}_{\left(  m,m-1\right)  }^{4}\right)  =\dim\widehat{Z}_{2t+2}%
^{2}\left(  \frak{g}_{\left(  m,m-1\right)  }^{4},\frak{g}_{\left(
m,m-1\right)  }^{4}\right)  =m-t-2
\]
In particular
\[
\dim F_{0}\widehat{Z}^{2}\left(  \frak{g}_{\left(  m,m-1\right)  }%
^{4},\frak{g}_{\left(  m,m-1\right)  }^{4}\right)  =m^{2}-2m+1
\]
\end{corollary}

\begin{proposition}
For any $m\geq4,$ the cohomology class $\left[  \psi_{k,r}\right]  $ of the
cocycles $\psi_{k,r}$ $\left(  1\leq r\leq2m-4\right)  $ \ are nontrivial and
linearly independent.\newline Moreover, $F_{0}\widehat{H}^{2}\left(
\frak{g}_{\left(  m,m-1\right)  }^{4},\frak{g}_{\left(  m,m-1\right)  }%
^{4}\right)  $ is a $\left(  m^{2}-3m+2\right)  $-dimensional \ subspace of
$F_{0}H^{2}\left(  \frak{g}_{\left(  m,m-1\right)  }^{4},\frak{g}_{\left(
m,m-1\right)  }^{4}\right)  $.
\end{proposition}

\begin{proof}
From the basis obtained in the previous proposition it follows easily that
$F_{0}\widehat{B}^{2}\left(  \frak{g}_{\left(  m,m-1\right)  }^{4}%
,\frak{g}_{\left(  m,m-1\right)  }^{4}\right)  =\left\{  \psi_{k,0}\;|\;2\leq
k\leq m\right\}  $. Observe in particular that most coboundaries have been
previously excluded from the cohomology space $F_{0}\widehat{Z}^{2}\left(
\frak{g}_{\left(  m,m-1\right)  }^{4},\frak{g}_{\left(  m,m-1\right)  }%
^{4}\right)  $ by property $\left(  P\right)  $.
\end{proof}

\begin{corollary}
For any $1\leq r\leq2m-4$ \ the cocycles $\psi_{k,r}$ $\left(  2\leq
k\leq\left[  \frac{2m-r}{2}\right]  \right)  $ form a basis of $\widehat
{H}^{2}\left(  \frak{g}_{\left(  m,m-1\right)  }^{4},\frak{g}_{\left(
m,m-1\right)  }^{4}\right)  $.
\end{corollary}

\begin{proof}
Follows from the previous result, as for $r\geq1$
\[
\widehat{H}_{r}^{2}\left(  \frak{g}_{\left(  m,m-1\right)  }^{4}%
,\frak{g}_{\left(  m,m-1\right)  }^{4}\right)  =\widehat{Z}_{r}^{2}\left(
\frak{g}_{\left(  m,m-1\right)  }^{4},\frak{g}_{\left(  m,m-1\right)  }%
^{4}\right)
\]
\end{proof}

\begin{proposition}
For $k\geq1$, any coycle $\psi_{k}\in\widehat{H}_{k}^{2}\left(  \frak{g}%
_{\left(  m,m-1\right)  }^{4},\frak{g}_{\left(  m,m-1\right)  }^{4}\right)  $
which satisfies $\psi_{k}\left(  C^{1}\frak{g}_{\left(  m,m-1\right)  }%
^{4},C^{1}\frak{g}_{\left(  m,m-1\right)  }^{4}\right)  =0$ is linearly
expandable if and only if $\psi_{k}=\lambda\psi_{2,k},\;\lambda\in
C\backslash\{0\}$.\newline Moreover $\frak{g}_{\left(  m,m-1\right)  }%
^{4}+\lambda\psi_{2,k}\simeq\frak{g}_{\left(  m,m-1\right)  }^{4}+\psi_{2,k}$
for any $\lambda\neq0$.
\end{proposition}

\begin{theorem}
Let $\psi_{k}\in\widehat{H}_{k}^{2}\left(  \frak{g}_{\left(  m,m-1\right)
}^{4},\frak{g}_{\left(  m,m-1\right)  }^{4}\right)  $ by a linearly expandable
cocycle which satisfies $\psi_{k}\left(  C^{1}\frak{g}_{\left(  m,m-1\right)
}^{4},C^{1}\frak{g}_{\left(  m,m-1\right)  }^{4}\right)  =0$. Then the
deformation $\frak{g}_{\left(  m,m-1\right)  }^{4}+\psi_{k}$ \ is isomorphic
to the nilradical of a solvable, complete Lie algebra $\frak{r}_{\left(
m,m-1\right)  }^{4,k}$ of dimension $2m+2$. Moreover, if $k=2m-5$ or
$k=2m-4$,then the Lie algebra $\frak{r}_{\left(  m,m-1\right)  }^{4,k}$ is rigid.
\end{theorem}

\begin{proof}
We have seen that the rank of the algebra $\frak{g}_{\left(  m,m-1\right)
}^{4}$ is $2$ for any $m\geq4$. From the previous proposition we know that
$\psi_{k}=\lambda\psi_{2,k}$ for $\lambda\neq0$. Without loss of generality we
can take $\lambda=1$. Let $f\in Der\left(  \frak{g}_{\left(  m,m-1\right)
}^{4}+\psi_{2,k}\right)  $ such that
\[
f\left(  X_{t}\right)  =\sum_{j=1}^{2m+2}f_{t}^{j}X_{j},\;1\leq t\leq2m+2
\]
The cocycle $\psi_{2,k}$ implies in particular
\[
\psi_{2,k}\left(  X_{2},X_{j}\right)  =X_{1+j+k},\;3\leq j\leq2m-1-k
\]
Taking the relation $\left[  f\left(  X_{2}\right)  ,X_{3}\right]  +\left[
X_{2},f\left(  X_{3}\right)  \right]  =f\left(  X_{4+k}\right)  $, we obtain
\[
f_{2}^{2}=\left(  1+k\right)  f_{1}^{1}%
\]
This proves that the deformation $\frak{g}_{\left(  m,m-1\right)  }^{4}%
+\psi_{2,k}$ has rank one. Now let $\frak{r}_{\left(  m,m-1\right)  }^{4,k}$
be the solvable Lie algebra whose Cartan-Maurer equations on the basis
$\left\{  \omega_{1},..,\omega_{2m+1},\theta\right\}  \;$are
\begin{align*}
d\omega_{1}  &  =\theta\wedge\omega_{1}\\
d\omega_{2}  &  =\left(  k+1\right)  \theta\wedge\omega_{2}\\
d\omega_{j}  &  =\omega_{1}\wedge\omega_{j-1}+\left(  k+j-1\right)
\theta\wedge\omega_{j},\;3\leq j\leq3+k\\
d\omega_{4+k+j}  &  =\omega_{1}\wedge\omega_{3+k+j}+\omega_{2}\wedge
\omega_{3+j}+\left(  k+5+j\right)  \theta\wedge\omega_{4+k+j},\;0\leq
j\leq2m-4-k\\
d\omega_{2m+1}  &  =\sum_{j=2}^{m}\left(  -1\right)  \omega_{j}\wedge
\omega_{2m+1-j}+\left(  2k+2m-1\right)  \theta\wedge\omega_{2m+1}\\
d\theta &  =0
\end{align*}
It is immediate to verify that this Lie algebra is solvable and complete. If
$k=2m-5$ or $2m-4$ $\ $the algebra is moreover ( see [1] for a standard proof
\ ). As we have $\frak{r}_{\left(  m,m-1\right)  }^{4,k}\simeq\left(
\frak{g}_{\left(  m,m-1\right)  }^{4}+\psi_{2,k}\right)  \oplus\frak{t}$, its
nilradical is clearly isomorphic to $\frak{g}_{\left(  m,m-1\right)  }%
^{4}+\psi_{2,k}$.
\end{proof}

\begin{remark}
The previous proposition provides us with a correspondence ( $k=2m-4,\;k=2m-5$
)%
\[
I:\widehat{H}^{2}\left(  \frak{g}_{\left(  m,m-1\right)  }^{4},\frak{g}%
_{\left(  m,m-1\right)  }^{4}\right)  \rightarrow\mathbb{Z}^{2m+1}%
\]
between the cohomology class of the cocycle $\psi_{2,k}$ and the $\left(
2m+1\right)  $-uples \newline $\left(  1,k+1,...,2k+2m-1\right)  ,$ where this
sequence gives the weight distribution of the eigenvalues of the action of the
torus $\frak{t}$ over $\frak{g}_{\left(  m,m-1\right)  }^{4}+\psi_{2,k}$.
\end{remark}

Let $\psi\in F_{1}\widehat{H}^{2}\left(  \frak{g}_{\left(  m,m-1\right)  }%
^{4},\frak{g}_{\left(  m,m-1\right)  }^{4}\right)  $ be a cocycle. Then this
cocycle is expressible as
\[
\psi=\psi_{1}+..+\psi_{r},\;\;\psi_{i}\in\widehat{H}_{i}^{2}\left(
\frak{g}_{\left(  m,m-1\right)  }^{4},\frak{g}_{\left(  m,m-1\right)  }%
^{4}\right)
\]

If $r$ is the greatest integer such that $\ \psi_{r}\neq0,$ then we call
$\psi_{r}$ the sill cocycle of $\psi$.

\begin{lemma}
For any cocycle $\psi\in F_{1}\widehat{H}^{2}\left(  \frak{g}_{\left(
m,m-1\right)  }^{4},\frak{g}_{\left(  m,m-1\right)  }^{4}\right)  $ the sill
cocycle $\psi_{r}\newline \left(  1\leq r\leq2m-4\right)  $ is linearly expandable.
\end{lemma}

\begin{remark}
This fact is already known for the deformations of the filiform Lie algebra
$L_{n}$ [21]. For this model the property is preserved, as it is a central
extension of the last algebra. The Lie algebra $\frak{g}_{\left(
m,m-1\right)  }^{4}+\psi_{r}$, where $\psi_{r}$ is the sill cocycle, is
usually called sill algebra.
\end{remark}

\begin{theorem}
Let $\psi\in F_{1}\widehat{H}^{2}\left(  \frak{g}_{\left(  m,m-1\right)  }%
^{4},\frak{g}_{\left(  m,m-1\right)  }^{4}\right)  $ be a linearly expandable
cocycle such that the sill cocycle $\psi_{r}$ satisfies \newline
1)\ $r\leq2m-5$\newline 2)\ $\psi_{r}\left(  C^{1}\frak{g}_{\left(
m,m-1\right)  }^{4},C^{1}\frak{g}_{\left(  m,m-1\right)  }^{4}\right)
=0.$\newline If there exists a component $\psi_{j_{0}}$ of $\psi$ such that
$j_{0}\leq r-1$ and satisfying $\psi_{j_{0}}\left(  C^{1}\frak{g}_{\left(
m,m-1\right)  }^{4},C^{1}\frak{g}_{\left(  m,m-1\right)  }^{4}\right)  \neq0$,
then the Lie algebra $\frak{g}_{\left(  m,m-1\right)  }^{4}+\psi$ is
characteristically nilpotent.
\end{theorem}

\begin{proof}
As the sill cocycle is itself linearly expandable, we know from previous
results that the algebra $\left(  \frak{g}_{\left(  m,m-1\right)  }^{4}%
+\psi_{r}\right)  $ has rank one. If $f\in Der\left(  \frak{g}_{\left(
m,m-1\right)  }^{4}+\psi\right)  $, it is easy to verify that the entries of
the matrix $\left(  f_{j}^{i}\right)  _{1\leq i,j\leq2m+1}$ of $f$ satisfy the
relations
\begin{align*}
f_{j}^{i}  &  =0,\;2\leq i<j\leq2m+1\\
f_{i}^{i}  &  \in\left\langle f_{1}^{1},f_{2}^{2}\right\rangle _{\mathbb{C}%
},\;3\leq j\leq2m+1
\end{align*}
Thus the sill cocycle implies%
\[
f_{2}^{2}=\left(  r+1\right)  f_{1}^{1}%
\]
Now let $\psi_{j_{0}}$ the summand of $\psi$ that satisfies $\psi_{j_{0}%
}\left(  C^{1}\frak{g}_{\left(  m,m-1\right)  }^{4},C^{1}\frak{g}_{\left(
m,m-1\right)  }^{4}\right)  \neq0$. Then there exists an $s\geq4$ such that
\[
\psi_{j_{0}}\left(  X_{3},X_{s}\right)  =\lambda_{3s}X_{1+s+j_{0}}%
\]
From this equation we obtain%
\[
\left(  s-1\right)  f_{1}^{1}+2f_{2}^{2}=\left(  1+s+j_{0}-2\right)  f_{1}%
^{1}+f_{2}^{2}%
\]
so that
\[
f_{2}^{2}=j_{0}f_{1}^{1}%
\]
and as $j_{0}\leq r-1$ we have $f_{1}^{1}=f_{2}^{2}=0$. It is not difficult to
see that the other summands of $\psi$ do not affect these relations. It
follows that the matrix of $f$ is strictly upper triangular, thus the
derivation $f$ is nilpotent and $\left(  \frak{g}_{\left(  m,m-1\right)  }%
^{4}+\psi\right)  $ is characteristically nilpotent.
\end{proof}

\begin{corollary}
For any $m\geq4$ and $k\geq3,\;1\leq s\leq r-1\leq2m-6$ the Lie algebra
\[
\frak{g}_{\left(  m,m-1\right)  }^{4}+\alpha\psi_{2,r}+\beta\psi
_{k,s},\;\;\alpha,\beta\in\mathbb{C-\{}0\}
\]
is characteristically nilpotent.
\end{corollary}

\section{Extensions of $\frak{g}_{\left(  m,m-1\right)  }^{4}$}

In this last section we study how the previous deformations are compatible
with central extensions of $\frak{g}_{\left(  m,m-1\right)  }^{4}$ by
$\mathbb{C}$. We will obtain, for any $m\geq4$, characteristically nilpotent
Lie algebras of characteristic sequence $\left(  2m-1,2,1\right)  $ in
dimension $2m+2$.

For $m\geq4$ let us define the Lie algebras $\frak{g}_{\left(  m,m-1\right)
}^{4,1}$ whose Cartan-Maurer equations are
\begin{align*}
d\omega_{1}  &  =d\omega_{2}=0\\
d\omega_{j}  &  =\omega_{1}\wedge\omega_{j-1},\;3\leq j\leq2m\\
d\omega_{2m+1}  &  =\sum_{j=2}^{m}\left(  -1\right)  ^{j}\omega_{j}%
\wedge\omega_{2m+1-j}\\
d\omega_{2m+2}  &  =\omega_{1}\wedge\omega_{2m+1}+\sum_{j=2}^{m}\left(
-1\right)  ^{j}\left(  m+1-j\right)  \omega_{j}\wedge\omega_{2m+2-j}%
\end{align*}

\begin{remark}
It is clear that these algebras have characteristic sequence $\left(
2m-1,2,1\right)  $ for any $m\geq4$ and that they are $\left(  m-1\right)
$-abelian. The only property which got lost by the extension is the natural graduation.
\end{remark}

Moreover, let $\overset{-}{\psi}_{2,k}$ be the prolongation by zeros of
$\psi_{2,k}$, i.e.,
\[
\overset{-}{\psi}_{2,k}=\left\{
\begin{array}
[c]{cc}%
\psi_{2,k}\left(  X_{l},X_{m}\right)  & \text{if \ }1\leq l,m\leq2m+1\\
0 & \text{if \ }l=2m+2\text{ or }m=2m+2
\end{array}
\right.
\]

\begin{proposition}
The cocycle $\overset{-}{\psi}_{2,k}\in H^{2}\left(  \frak{g}_{\left(
m,m-1\right)  }^{4,1},\frak{g}_{\left(  m,m-1\right)  }^{4,1}\right)  $ is
linearly expandable if and only if $k=2m-5$ or $k=2m-4$.
\end{proposition}

\begin{proof}
The only manner in which the differential form $d\omega_{2m+2}$ is closed is
that the cocycle $\overset{-}{\psi}_{2,k}$ only affects the differential form
$d\omega_{2m}$ adding the exterior product $\omega_{2}\wedge\omega_{3}$ or
changing the differential forms $d\omega_{2m-1},d\omega_{2m}$ by adding,
respectively, the exterior products $\omega_{2}\wedge\omega_{3}$ and
$\omega_{2}\wedge\omega_{4}$. This corresponds to the indexes $k=2m-4$ and
$k=2m-5$.
\end{proof}

\begin{remark}
Observe that the prolongation by zeros of a cocycle will be linearly
expandable in the extension if it is compatible with the adjoined differential
form $d\omega_{2m+2}$. This result can be easily generalized for prolongations
of cocycles $\psi_{t,k}$ for arbitrary $t\geq3$, but we are only considering
the cocicles for which $t=2$.
\end{remark}

.

The two next theorems are direct consecuences of the results of the previous
sections :

\begin{theorem}
For $m\geq4$ let $\frak{r}_{\left(  m,m-1\right)  }^{4,1,2m-5}$ be the Lie
algebra whose Cartan-Maurer equations are
\begin{align*}
d\omega_{1}  &  =\theta\wedge\omega_{1}\\
d\omega_{2}  &  =\left(  2m-4\right)  \theta\wedge\omega_{2}\\
d\omega_{j}  &  =\omega_{1}\wedge\omega_{j-1}+\left(  2m-6+j\right)
\theta\wedge\omega_{j},\;3\leq j\leq2m-2\\
d\omega_{2m-1}  &  =\omega_{1}\wedge\omega_{2m-2}+\omega_{2}\wedge\omega
_{3}+\left(  4m-7\right)  \theta\wedge\omega_{2m-1}\\
d\omega_{2m}  &  =\omega_{1}\wedge\omega_{2m-1}+\omega_{2}\wedge\omega
_{4}+\left(  4m-6\right)  \theta\wedge\omega_{2m}\\
d\omega_{2m+1}  &  =\sum_{j=2}^{m}\left(  -1\right)  ^{j}\omega_{j}%
\wedge\omega_{2m+1-j}+\left(  6m-11\right)  \theta\wedge\omega_{2m+1}\\
d\omega_{2m+2}  &  =\omega_{1}\wedge\omega_{2m+1}+\sum_{j=2}^{m}\left(
-1\right)  ^{j}\left(  m+1-j\right)  \omega_{j}\wedge\omega_{2m+2-j}+\left(
6m-10\right)  \theta\wedge\omega_{2m+2}\\
d\theta &  =0
\end{align*}
This algebra is solvable, decomposable and rigid. Moreover, its nilradical is
isomorphic to $\frak{g}_{\left(  m,m-1\right)  }^{4,1}+\overset{-}{\psi
}_{2,2m-5}$.
\end{theorem}

\begin{theorem}
For $m\geq4$ the Lie algebra $e_{1}\left(  \frak{g}_{\left(  m,m-1\right)
}^{4}+\psi_{2,2m-4}\right)  $ defined by
\begin{align*}
d\omega_{1}  &  =d\omega_{2}=0\\
d\omega_{j}  &  =\omega_{1}\wedge\omega_{j-1},\;3\leq j\leq2m-1\\
d\omega_{2m}  &  =\omega_{1}\wedge\omega_{2m-1}+\omega_{2}\wedge\omega_{3}\\
d\omega_{2m+1}  &  =\sum_{j=2}^{m}\left(  -1\right)  ^{j}\omega_{j}%
\wedge\omega_{2m+1-j}\\
d\omega_{2m+2}  &  =\omega_{1}\wedge\omega_{2m+1}+\sum_{j=2}^{m}\left(
-1\right)  ^{j}\left(  m+1-j\right)  \omega_{j}\wedge\omega_{2m+2-j}%
\end{align*}
is characteristically nilpotent. Moreover, its characteristic sequence is
$\left(  2m-1,2,1\right)  .$
\end{theorem}


\begin{thebibliography}{99}
\bibitem{}J. M. Ancochea. \textit{On the rigidity of solvable Lie algebras,
}ASI Nato Serie C247, 403-445 (1986).

\bibitem {}J. M. Ancochea, O. R. Campoamor\textit{. On (n-5)- 2-abelian Lie
algebras, } Comm. Algebra, \textit{to appear}.

\bibitem {}J. M. Ancochea, O. R. Campoamor\textit{. Characteristically
nilpotent Lie algebras, \ }Contribuciones matem\'{a}ticas n.3.

\bibitem {}J. M. Ancochea, O. R. Campoamor\textit{. Characteristically
nilpotent deformations of naturally graded Lie algebras of type }$Q_{n}$, preprint.

\bibitem {}F. Bratzlavsky. \textit{Sur les alg\`{e}bres admettant un tore des
d\'{e}rivations donn\'{e}, }J. of Algebra 30 (1974), 305-316.

\bibitem {}R. Carles. \textit{Sur le structure des alg\`{e}bres de Lie
rigides, }Ann. Inst. Fourier 34(3), 65-82 (1984).

\bibitem {}R. Carles. \textit{Sur les alg\`{e}bres caract\'{e}ristiquement
nilpotentes. }Publ. Univ. Poitiers 1984.

\bibitem {}C. Y. Chao. \textit{Uncountably many non isomorphic nilpotent Lie
algebras, }Proc. Amer. Math. Soc. 13 (1962), 903-906.

\bibitem {}J. Dixmier, W. G. Lister. \textit{Derivations of nilpotent Lie
algebras, }Proc. Amer. Math. Soc. 8 (1957), 155-157.

\bibitem {}J. Dixmier. \textit{Cohomologie des alg\`{e}bres de Lie
nilpotentes, }Acta Sci. Math. Szeged 16 (1955), 246-250.

\bibitem {}M. Gerstenhaber. \textit{On the deformations of rings and algebras.
}Ann. Math. (2) 79 (1964), 59-103.

\bibitem {}M. Goze. \textit{Crit\'{e}res cohomologiques pour la rigidit\'{e}
de lois algebriques, }Bull. Soc. Math. Belgique 43 (1991), 33-42.

\bibitem {}M. Goze, Yu. B. Khakimdjanov. \textit{Nilpotent Lie algebras,
}Kluwer Ac. Press 1996.

\bibitem {}N. Jacobson. \textit{A note on automorphisms and derivations of Lie
algebras, }\ Proc.Amer. Math. Soc. 6 (1955), 281-283.

\bibitem {}Yu. B. Khakimdjanov. \textit{Vari\'{e}t\'{e} des lois
d'alg\`{e}bres de Lie nilpotentes, }Geometriae Dedicata 40 (1991), 269-295.

\bibitem {}Yu. B. Khakimdjanov. \textit{Characteristically nilpotent Lie
algebras, }Math. USSR Sbornik 70 (1990), n.1.

\bibitem {}Yu. B. Khakimdjanov. \textit{On characteristically nilpotent Lie
algebras, }Soviet Math. Dokl. 41 (1990), n.1.

\bibitem {}G. F. Leger. \textit{A note on the derivations of Lie algebras.
}Proc. Amer. Math. Soc. 4 (1953), 511-514.

\bibitem {}G. F. Leger, S. T\^{o}g\^{o}. \textit{Characteristically nilpotent
Lie algebras, }Duke Math. J. 26 (1959), 623-628.

\bibitem {}F. J. Murray. \textit{Perturbation theory and Lie algebras, }J.
Math. Phys. 3 (1962), 451-468.

\bibitem {}A. Nijenhuis, R. W. Richardson. \textit{Deformations of Lie algebra
structures, }J. Math. Mech. 17 (1967), 89-105.

\bibitem {}A. Nijenhuis, R. W. Richardson. \textit{Cohomology and deformations
in graded Lie algebras, }Bull. Amer. Math. Soc. 72 (1966), 1-29.

\bibitem {}T. Skjelbred, T. Sund. \textit{Sur la classification des
alg\`{e}bres de Lie nilpotentes, }C. R. Acad. Sci. Paris 286 (1978), 241-242.

\bibitem {}S. T\^{o}g\^{o}. \textit{Outer derivations od Lie algebras, }Trans.
Amer. Math. Soc. 128 (1967), 264-276.

\bibitem {}S. T\^{o}g\^{o}. \textit{On the derivation algebras of Lie
algebras, }Cand. J. Math. 13(2), 201-216 (1961).

\bibitem {}M. Vergne. \textit{Vari\'{e}t\'{e} des alg\`{e}bres de Lie
nilpotentes, }These 3$^{eme}$ \textit{cycle, Paris 1966.}

\bibitem {}M. Vergne. \textit{Cohomologie des alg\`{e}bres de Lie nilpotentes.
Applications a l'\'{e}tude de la vari\'{e}t\'{e} des alg\`{e}bres de Lie
nilpotentes, }Bull. Soc. Math. France 98 (1970), 81-116.

\bibitem {}S. Yamaguchi. \textit{On some classes of nilpotent Lie algebras and
their automorphism group, }Mem. Fac. Sci. Kyushu Univ. Ser. A 34 (1981), 241-251.

\bibitem {}L. Zhu \textit{The construction of some solvable complete Lie
algebras, }J. Nanjing Univ. Math. Biquat. 15(1), 34-40 (1998).
\end{thebibliography}
\end{document}